# Balancing Accuracy and Complexity in Optimisation Models of Distributed Energy Systems and Microgrids with Optimal Power Flow: A Review

*Ishanki A. De Mel[a], Oleksiy V. Klymenko[a], Michael Short[a]\**

[a]Department of Chemical and Process Engineering, University of Surrey, Guildford, United Kingdom GU2 7XH
*m.short@surrey.ac.uk

**Abstract**

Optimisation and simulation models for the design and operation of grid-connected distributed energy systems (DES) often exclude the inherent nonlinearities related to power flow and generation and storage units, to maintain an accuracy-complexity balance. Such models may provide sub-optimal or even infeasible designs and dispatch schedules. In DES, optimal power flow (OPF) is often misrepresented and treated as a standalone problem. OPF consists of highly nonlinear and nonconvex constraints related to the underlying alternating current (AC) distribution network. This aspect of the optimisation problem has often been overlooked by researchers in the process systems and optimisation area. In this review we address the disparity between OPF and DES models, highlighting the importance of including elements of OPF in DES design and operational models to ensure that the design and operation of microgrids meet the requirements of the underlying electrical grid. By analysing foundational models for both DES and OPF, we identify detailed technical power flow constraints that have been typically represented using oversimplified linear approximations in DES models. We also identify a subset of models, labelled DES-OPF, which include these detailed constraints and use innovative optimisation approaches to solve them. Results of these studies suggest that achieving feasible solutions with high-fidelity models is more important than achieving globally optimal solutions using less-detailed DES models. Recommendations for future work include the need for more comparisons between high-fidelity models and models with linear approximations, and the use of simulation tools to validate high-fidelity DES-OPF models. The review is aimed at a multidisciplinary audience of researchers and stakeholders who are interested in modelling DES to support the development of more robust and accurate optimisation models for the future.

*Word count: 14,301*

**Key words:** Optimisation, Microgrids, Distributed Energy Systems, Nonlinear, Optimal Power Flow

**Nomenclature**

| | |
|---|---|
| AC | Alternating Current |
| ADMM | Alternating Direction Method of Multipliers |
| ADP | Approximate Dynamic Programming |
| CHP | Combined Heat and Power |
| BESS | Battery Energy Storage Systems |
| DC | Direct Current |
| DES | Distributed Energy Systems |
| DER | Distributed Energy Resources |
| DP | Dynamic Programming |
| DSM | Demand Side Management |
| GA | Genetic Algorithm |
| HV | High Voltage |
| LP | Linear Programming |
| LV | Low Voltage |
| MILP | Mixed-Integer Linear Programming |
| MINLP | Mixed-Integer Nonlinear Programming |
| MISOCP | Mixed-Integer Second Order Cone Programming |
| MIQP | Mixed-Integer Quadratic Programming |
| MPC | Model Predictive Control |
| MV | Medium Voltage |
| OPF | Optimal Power Flow |
| PSO | Particle Swarm Optimisation |
| SOCP | Second Order Cone Programming |
| SQP | Sequential Quadratic Programming |

## 1. Introduction

With increasing energy demands and pressure to reduce carbon emissions, distributed energy systems (DES) are predicted to play a vital role by 2050 in the energy industry [1]. These systems consist of small-scale distributed energy resources (DERs) located at or close to the premises of the end-user (known as the "prosumer" due to their ability to both produce and consume their own electricity). The term DES serves as an umbrella term for representative systems such as microgrids, energy hubs, distributed generation, smart local energy systems, and multi-carrier energy systems. Some systems, such as microgrids, can operate in both interconnected and isolated modes from external networks such as the electricity grid and gas networks [2]. There are many advantages to integrating these systems into existing energy networks, such as:

1) The localised nature of these systems makes it easier to integrate renewable energy resources and other generation technologies
2) DES technologies and their scale can be customised, depending on the local needs and energy availability
3) The ability to connect and disconnect from external networks can ensure uninterrupted power supply to critical loads such as hospitals [3].

DES have been increasingly investigated in the last two decades, predominantly through optimisation and simulation models [4]. Recent advances in this area include applications of Artificial Intelligence [5] and blockchain technologies [6], to enable trading between prosumers and external networks. However, in parallel to these advances, the fundamental components of these systems (such as generation/storage units used in DES, and connections to external networks) are being investigated in greater detail. DES components have inherent nonlinearities which are often ignored or linearised by modellers, such as AC power flow equations [7–12]. Results obtained from such models are often sub-optimal in comparison to the more-detailed model or even practically infeasible, meaning that they cannot be implemented. On the other hand, including these nonlinearities fully can produce intractable models.

DES optimisation models have generally been split into two categories, based on utility: design and operational models. Design models determine the optimal capacities of DERs and suitable locations for unit placement based on economic or environmental objectives. Operational aspects are usually incorporated in these design models to accurately determine the power required at different time points and hence, the capacities of the units. Operational models decide on the strategy of dispatching and operating the units available (unit sizes and locations are now fixed, as pre-determined by the design). These models are usually employed in day-ahead planning (offline), although real-time (online) operational models are increasingly being presented in literature [13]. Operational models are commonly referred to as the Energy Management System (EMS) of a microgrid. It can be described as an entity or controller which ensures optimal unit commitment and dispatch [14,15]. Most of the proposed models are discrete-time models, making decisions over a finite period or horizon.

The formulations of design and operational models are usually mixed-integer and nonlinear in nature, where integer variables are used to make decisions such as selecting which equipment to install and switching equipment on/off. Optimisation (mathematical programming) approaches have been widely used, with mixed-integer linear programming (MILP) being the most popular approach for both design and operation [16]. These models either ignore or linearise nonlinearities associated with DERs and/or networks and are therefore easier to formulate and solve using commercial solvers. Other commonly used deterministic optimisation methods include linear programming (LP) [17,18], mixed-integer nonlinear programming (MINLP) [19,20], sequential quadratic programming (SQP) [21,22], dynamic programming [23], etc. While linear programming methods can guarantee global optimality, nonlinear programming methods may not necessarily give globally optimum solutions, i.e. solutions may be stuck in local optima. The complexity of formulating and solving deterministic nonlinear programming problems have also resulted in a surge of metaheuristic and AI-based techniques applied to DES optimisation problems such as particle swarm optimisation (PSO) [24,25], genetic algorithms (GA) [26,27], artificial bee colony [28], grey wolf optimisation [29], etc. (see [30] and [31] for a summary and recent works, respectively). Note that heuristic and metaheuristic methods do not require any gradient information (i.e., can be applied to black box models or models with non-differentiable objective function or constraints) and can be simpler to formulate and apply. However, they may involve many parameters that require tuning and can cause complications due to parameter interactions [32]. Furthermore, they cannot guarantee optimality [33], can have long solution times, and their inherent stochasticity can result in different solutions from the same initial conditions.

DES are also simulated on various software packages and tools, which generally focus on detailed representation of units (related to generation, storage, and networks) [4,34]. Examples of such tools include HOMER [35], TRNSYS [36], and EnergyPlus [37]. They can be used either for design [38–42] (via gradient-free optimisation solvers or case studies) and/or operation [43–45]. The Distributed Energy Resources - Customer Adoption Model (DER-CAM) [46] is a popular MILP optimisation tool that has been used in many studies [47–49]. Although simulation and optimisation tools have been widely adopted, modellers may prefer to use mathematical programming and metaheuristic approaches where the modeller has full control and knowledge of the mathematical representations and techniques used.

Another important aspect of DES optimisation is optimal power flow (OPF). This consists of highly nonlinear and nonconvex constraints related to the underlying electrical network. Electrical grids around the world predominantly generate, transmit and distribute alternating current (AC) power. OPF plays a critical role in AC networks to calculate power losses in branches, compute node voltages and angles in multiple phases, and ensure that the power flows respect physical limitations of the network, while meeting consumer demand under objectives such as minimum generation cost or minimum power losses. Many studies that investigate the design and operation of grid-connected microgrids and DES treat OPF as a standalone problem or tool for optimisation and control of power flows [11,50]. However, it has been increasingly argued that elements of OPF must be included in the design and operational models of DES (such as within the EMS) to ensure that solutions of these models do not violate constraints associated with the distribution network, and therefore become infeasible [51]. This aspect of the optimisation problem has often been overlooked by researchers in the process systems and optimisation area. Despite ongoing research and evidence (such as [52] and [53]) demonstrating that DC approximations are not suitable in place of AC power flow equations, they are still used in grid-connected DES design and operation models. On the other hand, researchers who focused on OPF in a microgrid have often neglected the unit commitment and dispatch aspects of the microgrid, using less-detailed representations of components such as distributed generators and storage. This has led to a disparity between OPF models and DES optimisation models. Furthermore, widely used DES simulation tools have potentially overlooked this aspect as well.

There is an increasing interest in addressing this disparity, evident from recent studies which are analysed in this literature review, to consolidate the gap between DES and OPF. Figure 1 below shows the trend observed in literature found on the Scopus database for design and operation of DES while considering power flow and grid constraints, using the following search:

*TITLE-ABS-KEY(("distributed energy system*" OR "distributed energy resource*" OR "distributed generation" OR microgrid* OR "energy hub" OR "energy management system") AND ("power flow" OR "grid constraint*" OR "network topology") AND (optimi\*ation OR simulation*) AND (design OR operation* OR planning OR scheduling)) AND (model) AND PUBYEAR > 1999 AND PUBYEAR < 2021 AND (DOCTYPE(ar) OR DOCTYPE(cp) OR DOCTYPE(re))*

The rapid growth in the number of documents observed from 2015 to 2019 is indicative not only of the rising interest in DES optimisation, but also in this problem.

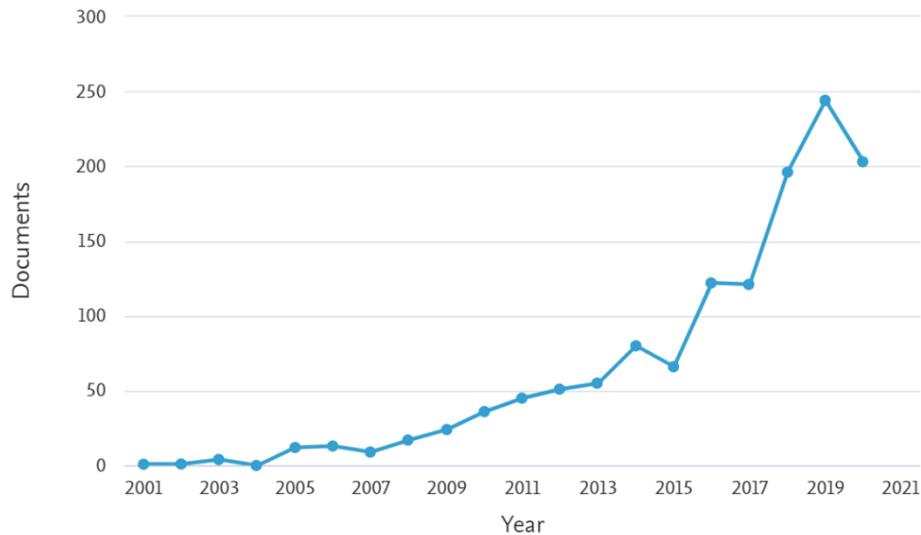

**Figure 1.** Trend of DES literature (2000 – 2020) containing targeted keywords such as "distributed energy system", "microgrid", "optimisation", "simulation", "power flow", and "network topology", as seen on Scopus.

*1.1    Novelty, research questions and scope*

While previously published reviews investigate and summarise either OPF/power flow approaches in smart grids and microgrids [50,54], or general DES/microgrid modelling approaches [30,55], this review differs from previous works by identifying and critically analysing a novel and promising group of studies that aim to:

- bridge the gap between OPF and DES design and operation
- solve nonlinearities and nonconvexities within a reasonable time frame with new approaches and techniques for the combined DES-OPF problem
- provide recommendations to future work on how to consolidate OPF and DES optimisation.

Throughout the review, we attempt to answer three questions:

1) What are the different levels of approximation evident through literature when modelling grid constraints in a DES or microgrid?
2) Should detailed power flow constraints be included in DES models?
3) Where applicable, how is OPF linked within DES models to produce optimal solutions within a reasonable time frame?

A more subjective research question is also explored, which may be of interest to DES modellers:

> Out of the modelling approaches explored, what level or method is most suitable to achieve a good complexity-accuracy balance for DES - OPF?

It is important to note that DES literature is extensive and there are many nonlinear DES models presented. We have restricted our overall search to studies from 2010 – 2021, paying more attention to recent work (2016 – 2021). We focus specifically on nonlinearities in the design stage and the EMS of a DES (in the operational stage) and ignore studies that do not include these aspects. While models for DC microgrids are abundant in literature, we focus on AC and AC-DC hybrid microgrids, considering existing grid infrastructure and ease of implementation [56]. Therefore, we focus on OPF in AC networks, while briefly discussing DC approximations in Section 2. We do not include literature concerning multi-microgrids, as this introduces other complexities due to interactions that may not be relevant to a single microgrid or DES (see [57] for more information). However, some of the literature presented here will be valuable to those considering multi-microgrid design and operation.

The novelty of this review is further emphasised in Figure 2. Firstly, we provide an overview of some fundamental models of DES (especially those that have laid foundations to more detailed modelling subsequently), and new contributions that do not consider detailed power flow constraints. These are labelled as *DES Baseline models*. We discuss reasons why modellers may choose to exclude such constraints, and potential repercussions. Next, we introduce the fundamentals of OPF and present a general formulation as seen in literature, including nonconvex functions for power flow. Finally, we investigate models that bridge the gap between OPF and DES design and operation and include detailed and transparent formulations. These are labelled as *Detailed and Combined (DES-OPF) models.* Some of these models tend to use novel and/or unconventional approaches to modelling. Figure 2 illustrates that the focus of this review is on the detailed and combined models which consolidate the disparity between DES design and operational models, and OPF models.

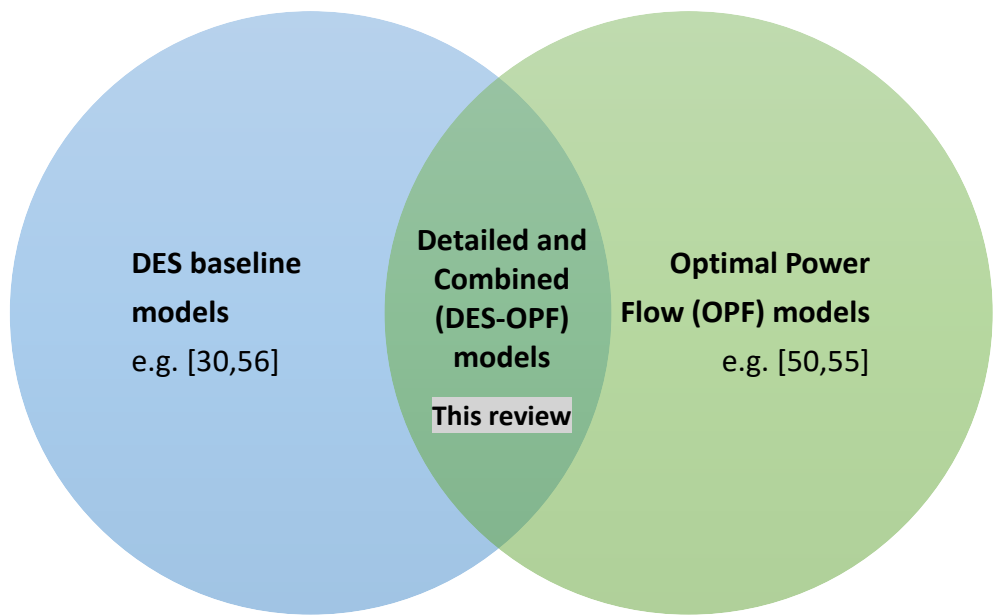

Figure 2. A graphical representation emphasising the novelty of this review, compared with existing reviews (references given as examples). Detailed and combined models are those that consider and consolidate both DES and OPF aspects.

Ultimately, this review should shed light to future modellers on considerations previously overlooked, and lead to more representative models suitable for real-time operation.

## 2. DES baseline models

The first category we consider are what we shall refer to as "DES baseline models". Most MILP DES models from the early 2010s fall into this category, although some of these models have undeniably laid the foundation for more complex models found more recently. This section aims to highlight some of these models, as well as selected recent studies that are of interest.

Ren and Gao [58] present an MILP model, predominantly for the design of a DES that allows the model to choose from a list of candidate DERs such as PVs, CHPs, fuel cells, etc., while considering operational constraints. The approach minimises total cost of the system, which can be broken down into capital costs, operational costs, carbon taxes, costs associated with trading with the grid, and so on. Note that this modelling framework has been used in and further developed by many subsequent publications. Linear constraints include the following:

- Supply and demand balances for electricity and heating loads
- Interactions with external networks, e.g. to prevent buying and selling electricity at the same time
- Capacity limits of technologies used
- Energy generation, storage, and use.

Nonlinearities associated with the operation of these technologies and the underlying distribution network are not included. Binary variables are used for assigning units, and operational states (such as start-up and shut-down). If there are bilinear terms where a product of binary and continuous variables is required, these are typically linearised using big-M formulations (as done in [59]), as linear solvers such as CPLEX [60] do not accept them. Chen et al. [61] consider a DES design model with renewable energy resources such as wind turbines and PVs, with particular focus on the sizing of energy storage systems within the DES. Although grid constraints are not considered in detail, the authors analyse how different storage capacities are chosen under two scenarios: 1) islanded operation, 2) grid-connected operation. This highlighted the need to consider the islanded operation mode within optimisation models, which had been overlooked in the past. Mehleri et al. [62] present an MILP model with novel constraints associated with a heating network. Each residential node is given the option of installing DERs and pipelines associated with the network. The study demonstrates the benefits of implementing a heating network to further reduce costs. Both Di Somma et al. [63] and Zhang et al. [64] evaluate the Pareto front of two objectives, economic and environmental, and how they influence the operational decisions of an MILP model. Silvente et al. [65] argue that model predictive control (MPC) strategies are useful in operational models in the presence of uncertain inputs (especially if there are intermittent renewable energy resources), allowing the model to react to any deviations from the initial conditions [66]. They present an MILP operational model that uses an MPC rolling horizon strategy to update data over short time intervals, based on renewable energy forecasts. Detailed grid constraints are not considered, despite enabling the microgrid studied to sell excess electricity to the grid. More recently, Karmellos and Mavrotas [67] have employed multi-objective optimisation and compared two different design approaches, one with predefined capacities for generation units and the other with only lower and upper capacity bounds. Constraints for a wide range of DERs are presented, however, they do not consider any inherent nonlinearities in these units and in the electrical network. In Huang et al. [68], DES design components are separated into layers of energy conversion, generation, and storage, based on graph theory. The overall model is solved using a two-stage MILP, where the first stage determines the optimal configuration of components and the second stage optimises the connections between the layers. The authors briefly discuss nonlinearities such as efficiencies of components but choose not to include them to avoid additional complexity.

Akter et al. [69] propose a distributed MILP model to incorporate peer-to-peer energy trading within the microgrid. The authors claim that the model is designed for an AC microgrid but do not present any constraints or nonlinearities pertaining to such networks.

To summarise, baseline DES design and operation formulations have been developed and improved over the years, by considering scenarios such as including two or more objectives, integrating heating networks, using MPC rolling horizon strategies, enabling peer-to-peer trading, considering grid-connected and islanded scenarios, etc. However, power flows and balances in all the models mentioned above remain over-simplified, and are commonly represented using the following equations or similar:

$$P_{i,t} = P_{i,t}^{Gen} - P_{i,t}^{Load} \quad (1)$$

$$Q_{i,t} = 0 \quad (2)$$

where $P_{i,t}$ and $Q_{i,t}$ represent the net active and reactive power, $P_{i,t}^{Gen}$ represents the generated power, and $P_{i,t}^{Load}$ represents the demand, each at a building or node $i$ with respect to time $t$. Linear approximations, specifically the DC approximation [70], has also been frequently employed instead of AC power flow constraints, as noted in the benchmarking models available in Sass et al. [71]. Considering node or bus $n \in \mathbf{N}$ and branch $(n,m) \in \mathbf{L}$, line susceptance $B_{nm}$ (the imaginary part of the admittance matrix) and voltage angles $\theta$ the power flow in branch $(n,m)$, denoted as $P_{nm,t}$, is calculated using Ohm's law, and the net power in node $n$, denoted as $P_{n,t}$, is calculated using Kirchhoff's current law:

$$P_{nm,t} = -B_{nm}(\theta_{n,t} - \theta_{m,t}) \quad (3)$$

$$P_{n,t} = \sum_{(n,m) \in \mathbf{L}} P_{nm} \quad (4)$$

. These equations assume that there is zero difference in voltage angles between branches (or buses), line resistances, and reactive power flows. Furthermore, voltages are assumed to be close to nominal voltage (which, in reality, can fluctuate, especially in the presence of intermittent renewable energy sources), and multiple phases and phase differences found in AC power flow are not taken into account. Essentially, only active power flows ($P_{n,t}$) are considered in DES baseline models. This is often unrealistic for AC power flow, as investigated by Purchala et al. [52], and can lead to infeasible results in the presence voltage fluctuations.

There are several reasons modellers may choose to exclude nonlinear and nonconvex grid constraints when modelling DES. Additional complexity is one of the key reasons, as nonlinearities are difficult to formulate and solve, with no guarantees on whether the solution is globally optimal or not. Linear approximations, such as the DC approximation, can also be employed in place of the nonlinear equations in the models to reduce the complexity. This applies not only to grid constraints, but also to other components such as batteries and CHPs, where aspects like battery ageing [72] and CHP efficiency [73] are sometimes best described using nonlinear representations. When considering many scenarios or samples for uncertainty or sensitivity analyses [74–76], it may be best to avoid such complexity as these models are already computationally expensive to solve. Furthermore, effective use of MPC strategies in day-ahead instances require fast solutions, for the optimisation problem is solved repeatedly over shifted time horizons. Additional complexity can hinder the process of obtaining results on time. It is also possible that modellers lack awareness of the impacts on the model outputs when some important nonlinearities are included. This may explain why detailed grid constraints have been largely ignored in many DES studies, despite the availability of relaxations and approximations that can potentially reduce the computational complexity.

## 3. Overview of OPF

OPF models are usually employed in transmission networks to minimise power losses within the network, or minimise costs for generation, while taking into account the topology and any physical limitations imposed by the network infrastructure [77]. When optimising, they do not typically consider constraints related to the various DERs, as done in DES design and operation models [78]. Rather, they consider the balances of total active and reactive power as functions of complex voltage at each node or branch, and any bounds associated with power, voltage, current, etc. Often, components such as distributed generators and energy storage are modelled in less detail [79]. With the increasing penetration of DERs and the implementation of microgrids, researchers have been investigating the use of OPF to determine power flows in distribution networks as well [80]. The distribution network is now described as 'active', i.e., there are some nodes which can now inject power, thus resulting in bidirectional power flow, rather than all nodes passively consuming power [79].

The reader is directed to Frank and Rebennack [77], which provides a comprehensive introduction to optimal power flow. To summarise, the two main nonlinear and nonconvex representations identified for single-phase networks are the bus injection model and branch flow model. There are many variants for each of these, depending on the coordinate system used (rectangular, polar, or combination) to define and compute the complex terms, relaxation techniques employed to convexify the overall formulation, and any additional constraints employed to calculate current magnitudes, represent components such as tap-changing transformers, etc. As the names of the two models suggest, the decision variables in the bus injection model are related to the nodes or buses, primarily nodal voltage and net active/reactive powers, while branch active/reactive powers and current form those of the branch flow model [81]. A basic OPF formulation based on the bus injection model, which is arguably the most commonly used formulation, is given below using Eqs. (5) - (14). Polar coordinates are used to represent voltage magnitude $V_{n,t}$, and angle $\theta_{n,t}$ over buses $n \in \mathbf{N}$, of which some have controllable generators $\mathbf{G} \subseteq \mathbf{N}$. The subscript $t$ has been included to represent discrete time, as modellers are often interested in net active, $P_{n,t}$, and reactive, $Q_{n,t}$, power flows over a time horizon. $G_{nm}$ and $B_{nm}$ represent line conductance and susceptance over branch $(n, m)$, respectively. Note that we do not consider multiple phases in the formulation below, but it is nevertheless an important aspect of OPF.

$$\min \sum_{n \in \mathbf{G}} C_n(P_{n,t}^{Gen}) \tag{5}$$

$$P_{n,t} = P_{n,t}^{Gen} - P_{n,t}^{Load} \quad \forall n \in \mathbf{N} \tag{6}$$

$$Q_{n,t} = Q_{n,t}^{Gen} - Q_{n,t}^{Load} \quad \forall n \in \mathbf{N} \tag{7}$$

$$P_{n,t} = V_{n,t} \sum_{m=1}^{N} V_{m,t}((G_{nm}\cos(\theta_{n,t} - \theta_{m,t}) + (B_{nm}\sin(\theta_{n,t} - \theta_{m,t})) \quad \forall n \in \mathbf{N} \tag{8}$$

$$Q_{n,t} = V_{n,t} \sum_{m=1}^{N} V_{m,t}((G_{nm} \sin(\theta_{n,t} - \theta_{m,t}) - (B_{nm} \cos(\theta_{n,t} - \theta_{m,t})) \quad \forall n \in \mathbf{N} \quad (9)$$

$$P_n^{Gen,min} \leq P_{n,t}^{Gen} \leq P_n^{Gen,max} \quad \forall n \in \mathbf{G} \quad (10)$$

$$Q_n^{Gen,min} \leq Q_{n,t}^{Gen} \leq Q_n^{Gen,max} \quad \forall n \in \mathbf{G} \quad (11)$$

$$V_n^{min} \leq V_{n,t} \leq V_n^{max} \quad \forall n \in \mathbf{N} \quad (12)$$

$$\theta_n^{min} \leq \theta_{n,t} \leq \theta_n^{max} \quad \forall n \in \mathbf{N} \quad (13)$$

$$V_{slack} = V_{nominal}, \qquad \theta_{slack} = 0 \quad (14)$$

The objective function in Eq. (5) minimises total cost as a function of power generated $P_{n,t}^{Gen}$. Eqs. (6) and (7) calculate the net active and reactive powers in each bus, respectively, where the superscript $Load$ indicates demand and $Gen$ indicates generated active and reactive power. Eqs. (8) and (9) are used to calculate the voltage and angle, as net active and reactive power are functions of them. Inequality constraints in Eqs. (10) - (13) ensure that active power, reactive power, voltage magnitude, and angle remain within predefined upper and lower bounds. Finally, Eq. (14) ensures that the voltage is fixed to nominal and angle remains zero at the system slack bus, referenced using the subscript $slack$. The slack bus serves as a reference point in the network and ensures a feasible solution as it has no limits on positive or negative power injections. Note that isolated microgrids, which are not connected to a larger/external grid, do not typically have a slack bus [82].

As evident from the above formulation, Eqs. (8) and (9) are highly nonconvex; their solutions are often only locally optimal, and global optimality cannot be guaranteed. OPF models are usually formulated as NLP problems and solved using interior-point methods [83]. Methods to obtain convex relaxations of these equations may also be used, where the feasible region of the optimisation problem is expanded to a convex space which includes the original nonconvex feasible region. These have become increasingly popular due to the tractability and possibility of obtaining globally optimal solutions [84]. Relaxation methods vary depending on the power flow model used (branch flow, bus injection, etc.), the coordinate system used (polar or rectangular), and the preferred method of solving the optimisation problem (Semidefinite programming (SDP), Second Order Cone Programming (SOCP), etc.) [84]. While both OPF models, bus injection and branch flow, have convex relaxations, those of the branch flow model are often favoured as they eliminate phase angles of voltages and currents [85]. The reader is directed to [50] and [84] for details on these convex relations, such as the use of McCormick envelopes for bilinear terms [86], as the focus of this review is not on OPF alone. Approximations and relaxations currently employed by DES-OPF models are discussed in Section 4.

Despite the common assumption that the underlying distribution network of the DES of microgrid is balanced, this is not the case in reality. While this simplifies the formulation, distribution networks are very rarely balanced across three phases. If this level of detail is to be considered in the design or operational model, multi-phase unbalanced OPF formulations should be used. The bus injection model is typically avoided in multi-phase scenarios, and modifications to the branch flow model are preferred due to their superior numerical stability [87]. Other methods include the Z-bus method [88], where the Bus Impedance Matrix, also known as the Z-bus, is calculated by inverting the Admittance matrix (Y-bus matrix), which is generally required in OPF calculations [89]. Modelling and simulating multi-phase unbalanced networks are discussed further in Section 4.

## 4. Detailed and combined (DES-OPF) models

Firstly, the studies that have focused on combining DES design and OPF are identified and analysed, followed by studies combining DES operation (scheduling) and OPF. We explore design and operational models separately as design models tend to have more decision variables when compared to operational models. As discussed in Section 1, the design model variables are the DER capacities and choices available, along with some operational variables such as amounts of power and heat produced, consumed, stored, sold, etc. Note that all design variables are treated as parameters in the operational models. Regardless of the formulation used in the DES design or operational model, whether linear or nonlinear, the inclusion of nonconvex power flow equations (as seen in Section 3) makes the overall problem NP-hard. Therefore, some of the discussed models have employed novel approaches to maintain the fidelity of AC OPF, while reducing complexity.

Some studies have also linked optimisation models to simulation tools. Their use is discussed in greater detail following the analysis of the DES-OPF models.

*4.1    Analysis – Design*

To introduce the reader to the DES design – OPF paradigm, some of the earlier works focused on "Distributed Generation Planning" [16] for the optimal placement of distributed generators, such as Moradi and Abedini [90] and Kaur et al. [91]. There are key differences between these models compared to more recent DES and microgrid design models. While both consider the placement of conventional and renewable generators using discrete variables, the focus in the former is on distribution networks. Network-related objectives such as power losses are used and the individual users in the network are not considered. In DES models, there is greater focus on individual users and their preferences within the network, by including more user-focused objectives such as minimising design costs and aspects like selling to the grid and demand response schemes. Nevertheless, both the above-mentioned studies emphasise the importance of considering the underlying distribution network and the inclusion of both active and reactive power balances at each bus, voltage limits, thermal limits (as a function of apparent power) and equations for calculating voltage stability at each bus. Due to the presence of complex voltages, nonconvex nonlinearities, and discrete variables, a combination of metaheuristics or exact methods have been employed to solve the models. For instance, GA and PSO are utilised in Moradi and Abedini [90], and SQP with branch-and-bound in Kaur et al. [91]. The former study outperforms the use of either GA or PSO alone, and the latter outperforms the use of PSO alone, both in terms of computational speed and obtaining better objective solutions. This highlights that hybrid metaheuristic or exact techniques are preferred over lone metaheuristic solvers.

On the other hand, Foster et al. [7] investigate whether greater complexity is required when using power flow equations in Distributed Generation Planning models. Three mathematical formulations are presented: 1) an MINLP complete power flow model, incorporating reactive power, power losses due to line resistance, and complex voltages in rectangular form, 2) an MILP model with DC power flow approximations (similar to that presented in Section 2 of this review), and 3) an MINLP model containing similar DC linear approximations as 2, but incorporating quadratic functions for line resistances (labelled as DCQL). The objective function minimises the power imported from external sources, such as the electrical grid. As solving highly nonconvex models cannot guarantee global optimality, the authors recommend finding the objective value for a relaxed version of the model to obtain a potential lower bound for the global optimum. This is done by removing all the constraints apart from two: the integrality constraint for the binary variables determining the placement of the DERs, and a constraint that limits the maximum capital expenditure on the DES design (i.e. the stakeholders' budget). This lower bound is then used to calculate a relative gap, used as a measure of the model performance. Interestingly, results show that the DCQL performs very similarly to the complete AC model with respect to the objective and relative gap, suggesting that not all equations in the AC-OPF model are required in DES design and operation. Despite the valuable contributions of this study which are also applicable to future DES Design – OPF models, the authors do not discuss if and how DES placement and sizing are different in the two models. This is important because, while the designs obtained from the DCQL model may be feasible under those constraints, they may not be feasible under the full AC constraints. This warrants further investigation and validation to confirm that DCQL approximations can be used in place of the full AC equations. The study, however, does confirm that the linear DC model employed in many studies is a poor representation of OPF when the DES is connected to the AC external grid.

On the DES Design – OPF paradigm, Morvaj et al. [8] has led the investigation on how the inclusion of detailed grid constraints in their models impacts the optimal design and operation of DES. They directly address and attempt to consolidate the optimisation of technologies in DES with detailed power flow equations and constraints in their models. Three models for both design and operation are presented, under economic and environmental objectives: 1) a GA-based design model is linked to an MILP operational model (or "energy hub") and MATPOWER [92], 2) an MILP design and operational model with linearised AC power flow equations, and 3) a GA design model linked to an operational MILP with linearised AC power flow equations and MATPOWER. Note that MATPOWER serves as a post-optimisation check using the Newton-Raphson method, and cannot directly influence the optimal design and operation, whereas the linearised equations included in Model 2 and 3 can. The DES aspect used in this study has linear/linearised representations for the DER technologies used, similar to those discussed in Section 2. Building on from equivalent nonconvex expressions for power flow presented in Section 3, the authors have also included a nonlinear expression to calculate the current ($I$) in each branch ($n, m$) with respect to each timestep ($t$):

$$I_{nm,t} = \frac{\Delta V_{nm}}{Z_{nm}} = \frac{V_n(cos\,\theta_n) + \mathbf{i}V_n(sin\,\theta_n) - V_m(cos\,\theta_m) - \mathbf{i}V_m(\sin\,\theta_m)}{R_{nm} + iX_{nm}} \quad (15)$$

where $R$, $X$, and $Z$ are resistance, reactance, and impedance of the line, respectively, while $V$ and $\theta$ represent voltage magnitude and angle, respectively. Note that this is an additional nonconvex equation with complex terms. Without resorting to DC power flow equations, their linear approximation for power flow equations ensures line conductance ($G$) and the difference in voltage between the branches are not set to zero. This allows the calculation of reactive power flow between branches ($n, m$), which is assumed zero in the DC approximation.

A nominal current value can be obtained by linearising Eq. (15), as detailed in Morvaj et al. [8]. Linearising the complex terms requires separating the real and imaginary parts of this equation, and considering only the magnitude of current, followed by piecewise linear approximations for quadratic terms. Mashayekh et al. [78], who also implement a linearised power flow model within the DER-CAM design environment, also adopt linear approximations. However, these are different to those in Morvaj et al. [8] as both real and imaginary parts are retained. While these linear approximations retain much more of the AC power flow characteristics compared to the DC approximations, the DES networks in both studies have been assumed to be balanced, thus considering only a single phase in the distribution network. This is a significant limitation when considering the wider applicability of these studies, as distribution networks are unbalanced in reality, and the best practice would be to consider all three phases. Yang and Li [93] highlight the importance of deciding whether to use balanced or unbalanced models, as results can vary significantly between the two. In a DES context, this can impact feasibility of the design obtained. Although the MILPs produce lowest costs, guarantee global optimality, and have low relative errors when power flow calculations are verified with nonlinear versions, the designs from these studies are not tested under the full AC power flow model to ensure that they are practically feasible. This is another limitation of these studies. Nevertheless, a key finding from Morvaj et al. [8] is that the designs produced by the three models are different. Models 1 and 3 may not produce globally optimal solutions, but the designs produced by these models are more likely to be practically feasible. This presents an accuracy-complexity dilemma to the modeller, and if DES design models are to be used in real-world applications, feasibility of the resulting design is more important that achieving global optimality.

While the previously discussed models have either linearised the power flow constraints and/or used metaheuristics to solve the MINLP formulation, Jordehi et al. [94] retain the MINLP formulation for OPF and solve it using the mathematical programming solver DICOPT [95]. The main focus of this study is to find optimal placement for EV battery swapping stations within a microgrid with other DERs such as PVs and geothermal units, while considering AC OPF constraints in a balanced medium voltage network. Although the results of bus voltages and angles are not presented in detail, the study identifies that the placement of the battery swapping station can indeed impact the power flow losses and bus voltages. This highlights the importance of including OPF in DES Design, and the need to discuss how these additional constraints can impact the overall optimisation. While the authors do not discuss this in the study, it can be assumed that solving this problem via DICOPT, a solver that uses an Outer Approximation/Equality-Relaxation algorithm, can be difficult as linearising the master problem and solving an MILP may not be a valid outer approximation in a nonconvex problem [96]. Inclusion of initialisations and extensions to the solvers used (such as the use of a feasibility pump) may enable more modellers to use mathematical programming approaches to obtain more reliable results, rather than relying on metaheuristic approaches.

Focusing on medium voltage (MV) networks, Sfikas et al. [21] considers the design of a microgrid under grid-connected and island modes. The objectives include minimising total cost of energy, which is the ratio of total annualised cost and total energy supplied, and total annual energy loss. While the standard bus injection model is utilised for balanced nonlinear OPF in a grid-connected mode, additional inequalities are presented to prevent active and reactive power exchange across the substation during island mode. This means that the slack bus is severely restricted in this mode. Unlike other design studies discussed above and in *DES baseline models*, no discrete variables are used in the formulations of DERs, particularly battery storage. The absence of discrete decisions allows the authors to use SQP to solve their model. However, this results in a more simplistic representation in the DES aspect, which may not be realistic despite the greater detail in the OPF equations and limits the accuracy of the results. Nevertheless, results highlight that taking the island mode into account demands higher dispatchable capacity to ensure demand is met, and small power imbalances do exist within the system. While the study emphasises the importance of taking both island and grid-connected modes into account, it is evident that achieving high fidelity with respect to OPF can result in compromised accuracy of DES design and operational aspects. This is a challenge most modellers investigating the DES-OPF paradigm face.

*4.2 Analysis – Operation*

In the earlier works of DES optimisation, Chen et al. (2013) [97] consider some constraints for both DC and AC current in the network within an EMS model for DES operation. This is an interesting development in DES/microgrid optimisation because the microgrid in their test case transmits/distributes power as DC, which however is converted to AC for end-user consumption. Transmission capacity limits and sinusoidal functions for voltage and current calculations are used for AC power flow, without sole reliance on linear DC approximations (as seen in *DES baseline models* discussed in Section 2). However, there are several limitations associated with this study with respect to the DES-OPF paradigm. Firstly, the microgrid is of residential-commercial scale; however, it is unclear why microgrid connections to the transmission network are considered, when this network usually transmits high-voltage power compared to residential-commercial level distribution networks. Secondly, it does not explore AC and DC sub-grids and bus configurations within the microgrid, as done in AC-DC hybrid microgrid models to minimise costs and losses associated with power conversion [98,99]. Thirdly, a detailed formulation for power flow, as presented by Eq. (5) - (14) in Section 3, is not considered. Finally, we assume that the nonlinear sinusoidal functions described have been linearised as the model is solved using CPLEX, a linear solver that cannot solve nonlinear terms. Any approximations or relaxations, and the impacts to the results in doing so are not detailed. Nevertheless, the paper highlights the importance of considering AC power flow equations in DES operation.

More recently, studies have included more detailed formulations for both DES technologies and power flow, and have considered the type of network the DES or microgrid is connected to. However, as previously discussed in Section 4.1, many studies focus on modelling balanced networks. On the DES Operation – OPF paradigm, these include Nemati et al. [9], Chen et al. [100], Shuai et al. [12], Zafarani et al. [101], Shi et al. [11], and Soares et al. [102], which are analysed first in this section. On the other hand, further advancements are seen in Thomas et al. [10], Jin et al. [103], and Morstyn et al. [104], where unbalanced three-phase power flow is considered.

Similar to Morvaj et al. [8] in Section 4.1, Nemati et al. [9] argue that previous models in literature do not pay attention to network constraints and inclusion of combined active/reactive power in DES models. Lack of validation of various optimisation methods is also highlighted. Aligning with the DES Operation – OPF paradigm, they present a day-ahead model, solved using two different methods, an MILP and a GA. The MILP, which they have linked to MATPOWER [92] to calculate network losses, is of particular interest in this review. Note that MATPOWER only captures balanced networks and does not have the functionality for unbalanced networks, to the best of the authors' knowledge. As MATPOWER cannot influence optimal decisions, the authors of [9] propose a multi-layer algorithm which involves first solving the MILP with linearised DER cost functions, subsequently solving the power flow calculations on MATPOWER and triggering an iterative process between the MILP and simulation to ensure that there are no violations. This appears to be less computationally efficient than integrating linearised power flow equations in the MILP and using MATPOWER as a verification tool, as done by Morvaj et al. [8]. However, this does ensure that the MILP designs are feasible under OPF equations for a balanced network. The study makes other valuable contributions by analysing scenarios with different objectives, such as maintaining grid stability, reducing carbon emissions, and so on. The MILP performs better than the GA in most scenarios, with the exception of highly constrained scenarios. This, once again, suggests that it may be better to use exact (mathematical programming) methods or hybrid methods, rather than relying on metaheuristics alone.

Chen et al. [100] aim to study both electricity and heating distribution networks in multi-carrier energy systems; but rather than using tools such as MATPOWER in Nemati et al. [9], they embed SOCP relaxations of the bus injection model in their formulation. They discourage modellers from solving nonconvex MINLPs as they may not be tractable and assert that SOCP relaxations (which are convex) are much easier to solve. Due to the presence of discrete decisions in the DES aspect, the overall problem is transformed to a Mixed Integer Second Order Cone Programming (MISOCP) model. The advances of this study lie in the optimisation approaches used, as they apply partial surrogate cuts [105] to decompose the MISCOP linear (MILP) and nonlinear (SOCP) substructures. The linear substructure consists of continuous linear and nonlinear variables associated with linear constraints and the MILP master problem solves this to obtain the lower bound of the overall problem. This retains the fidelity of the model (as it eliminates the need to obtain linear approximations for some of the nonlinear power flow constraints) while achieving relatively fast solutions. Solutions to the discrete variables are obtained through the MILP, but the lack of nonlinearity may result in infeasible discrete solutions. This is relaxed using nonnegative slack variables and integer cuts to prevent the algorithm from returning the same discrete solutions. The SOCP subproblem contains the nonlinear substructure, which is solved to obtain the upper bound of the overall problem. Ultimately, the algorithm is solved until it achieves solutions at or below the optimality gap specified by the modeller. A case study with no network constraints shows that operational costs can be severely underestimated if these are not included, further confirming the need to include them in DES operation. Testing this algorithm on a large-scale case study confirmed the ability of the model to solve quickly, despite having nonlinearities, showing that the model can potentially be used in multiple microgrid/energy hub scenarios. However, it is unknown whether the SOCP relaxations will be valid for a three-phase unbalanced scenario, highlighting that much work still needs to be done in this area.

Shuai et al. [12] also investigate microgrid scheduling, and confirm that DC formulations are not suitable in an EMS when it is connected to the AC distribution network, especially because of the "higher resistance-to-reactance" ratio in microgrid power lines. They also recognize that detailed modelling of components, such as batteries, can improve the accuracy of the model even though nonlinearities may be introduced. The authors attempt to consolidate all these aspects by proposing an MINLP that is solved using Approximate Dynamic Programming (ADP) where the original problem is decomposed to sequential subproblems. Power flow constraints detailed in the study are similar to those of the bus-injection model, and rather than using any convex relaxations as done in Chen et al. [100], nonconvexity is retained in the formulation. As ADP cannot achieve global optimality, its performance is compared with PSO (a metaheuristic solver) and Dynamic Programming (DP). ADP achieves faster approximate results close to the DP global optimum, surpassing the PSO in both time and quality. It confirms that deterministic optimisation techniques provide better solutions than any approximate or metaheuristic techniques, and convex relaxations or linear approximations can be used if speed is valued over accuracy.

Zafarani et al. [101] model the operation of a low-voltage energy hub with EVs and CHP in the presence of uncertainties, with the aim of minimising operational cost. A Robust Optimisation approach is used to obtain the optimal schedule for the worst-case values for the uncertain parameters. While considering the bus injection OPF equations, they also include nonlinear constraints for a natural gas network. Similar to Sfikas et al. [21] under DES Design – OPF, no discrete variables are included in the formulations, thus ignoring the on-off and ramp-up/down statuses of the CHP. Compared to other DES formulations, the exclusion of binary variables reduces the complexity of solving the problem; however, it also results in a less realistic representation. With respect to the nonconvex power flow constraints, piecewise and circular-plane linear approximations are used to formulate it as an LP. The globally optimal solutions of this are compared with an NLP that provides locally optimal but OPF-feasible solutions. Comparisons between the two show that the NLP requires more iterations (as expected), with LP calculation errors for the electrical network and related constraints at 2.5% or less. Unfortunately, the authors do not share if or how the operational schedule varies between the two deterministic models, as it is unlikely that the same results would be produced by both models. Overall, it is evident that, although the LP provides faster solutions for uncertainty analysis, the accuracy has been compromised to do so and therefore may not produce practically feasible schedules.

Continuing on studies that are based on low-voltage balanced networks, unlike the studies discussed above, Shi et al. [11] propose a model in which a main microgrid controller communicates with several local controllers that communicate with one another. The objective function aims to minimise operational costs (including those incurred because of dissatisfied customers) and power losses, incorporated as weights that can be set by the modeller. Instead of the bus injection model for OPF, the branch flow model is used where the following quadratic equality constraint is relaxed to be an inequality:

$$\left(I_{nm,t}\right)^2 = \frac{\left(P_{nm,t}\right)^2 + \left(Q_{nm,t}\right)^2}{\left(V_{n,t}\right)^2} \tag{16}$$

where the square terms $\left(I_{nm,t}\right)^2$ and $\left(V_{n,t}\right)^2$ are replaced by equivalent variables $i_{nm,t}$ and $v_{nm,t}$. This relaxation is claimed to be exact if the voltage and power imported from the main grid are within specified bounds and is verified by the authors for the results of the test case. To solve these constraints and the objective function, the authors use the predictor-corrector proximal multiplier (PCPM) algorithm, which is typically used to decompose and solve large-scale distributed OPF problems [106]. Run times for the total horizon considered are not reported; however, these algorithms can produce fast results. The downside is that they can only be used when the overall problem is convex and may not be suitable for analysis of three-phase unbalanced networks. Further studies are required to compare with the original nonconvex formulation and ensure that the convexification does not result in errors under more extreme conditions, such as in the presence of voltage fluctuations.

With respect to medium voltage (MV) networks, Soares et al. [102] model day-ahead operation in a MV balanced distribution network. Note that high voltages (HV) are atypical of distribution networks; they are known as transmission networks. Two objectives are considered, minimising operational cost, and maximising the power reserve (for customers with critical loads). HV/MV transformer limits and line thermal limits are also included, in addition to the bus injection model formulation presented in Section 3. The algorithm employed is a hybrid MILP-PSO, where the MILP is solved with DC power flow equations while the PSO retains the nonconvex AC power flow constraints. The optimal solution obtained from solving the MILP is fed into the PSO as an initial swarm solution to obtain feasible solutions. The authors utilise cluster computing to solve this algorithm within 30 mins, which would otherwise take a significantly longer time. When compared to an MILP with DC power flow, the authors claim that this model results in power losses and voltage violations, however, the relative errors are not reported. Nevertheless, it signifies that finding feasible solutions is valued over finding globally optimal solutions in this instance, which is an important consideration when formulating DES Operation – OPF models.

A pioneering study in the DES Operation – OPF paradigm has been conducted by Thomas et al. [10], who address the disparity between OPF and EMS models by integrating a Power Quality framework in their MILP operational model. While all the previously discussed studies consider a balanced network, this study uses OpenDSS [107] as a simulation tool within their algorithm to capture the voltage unbalance across phases in unbalanced distribution networks. The study also analyses the provision of additional services by the microgrid to the distribution network, such as providing active power to meet demand, provided that power quality is maintained. Interestingly, rather than using detailed power flow models within the formulation, indices for changes in voltage (kept within limits), harmonics (that may interfere with grid stability), and Voltage Unbalance Factor (VUF) (as the AC networks typically have three phases across which all loads are distributed) are calculated using OpenDSS and MATLAB [108]. Furthermore, to perform these calculations, the microgrid loads are divided into categories (such as lighting, appliances, etc.) and active and reactive power balances are described for each. Note that these authors also use power factors to calculate reactive power corresponding to loads and energy storage. The objective function, to minimise operational cost, also includes a penalty if reactive power is generated within the microgrid. If the solution generated by the MILP violates the limits for the indices, the MILP is updated with additional constraints to prevent violations in the subsequent iterations, thus ensuring feasibility of the operational schedules proposed by the model. These constraints appear nonlinear, so incorporating them requires linear approximations which are not discussed in the paper. The algorithm solves within seconds (even in a scenario with three iterations), showing immense potential for day-ahead modelling considering power flow calculations.

Similar to the study by Thomas et al. [10] above, Jin et al. [103] also demonstrates the use of OpenDSS for three-phase power flow calculations within a DES-OPF day-ahead framework. The authors also introduce discrete variables within the optimisation model to represent Remote Controlled Switches that can change the topology of the underlying distribution network. This is an interesting feature that has not been included in the previously discussed studies. A hybrid algorithm based on GA and NLP is used to solve the DES with a nonlinear natural gas network representation for heating. GA is used to solve subproblems that make discrete decisions, and the results of these are fed to the NLP, resulting in an iterative process until the convergence criteria are met. It is unclear how OpenDSS influences the optimisation algorithm in this study, and what criteria exist to detect violations. However, results show that when voltage violations are detected, network topology reconfiguration takes place to eliminate violations and enable normal operation. Thus, the inclusion of network switches can be very useful when networks can get congested during peak times, and changes in topology can ensure voltage and capacity limits are still being met.

Another important development in DES Operation – OPF paradigm is presented by Morstyn et al. [104]. The authors present OPEN, a Python-based open-source platform that provides a flexible modelling environment to users without the constraints of existing Energy Management tools. Both linear and nonlinear DES technology representations are supported, particularly for energy storage, facilitating their comparison to establish the accuracy-complexity balance required for operational models. One of the key achievements of this platform is the easy interlinking of OPF tools for multi-period DES optimisation with commonly used operational strategies such as MPC (as discussed in Section 2), giving the user the option of testing both balanced and unbalanced network scenarios. PandaPower [109], a user-friendly version of the python OPF tool PyPower (based on MATPOWER [92]), is used for static, balanced and nonlinear power flow calculations within DES models. As there are no three-phase unbalanced OPF tools available on Python, the authors introduce the *Network_3ph* class, containing nonlinear formulations for power flow simulation and linearised formulations for power flow optimisation. Rather than using either the bus injection model or the branch flow model as previously discussed studies have done, both the linear and nonlinear formulations are based on the Z-Bus method [88,89] (mentioned first in Section 3) and specifically suited for operational models. An advantage of the formulation employed for the *Network_3ph* class is that it allows three-phase distribution transformer connections, such as wye and delta, to be included in the Y-bus matrix without resulting in singularities [89]. Comparing the performance of the nonlinear *Network_3ph* formulation with the OpenDSS simulation tool [107] shows a high degree of accuracy, which confirms that this will be a valuable tool in the development of the DES operation – OPF paradigm. Note that the applicability of the *Network_3ph* to design models is unknown at present and requires further investigation. While a three-phase test case with electric vehicles and charging points has been used to demonstrate the use of this class, no comparisons have been performed against the linearised-nonlinear formulations, and balanced-unbalanced formulations. Such a study would motivate future modelling efforts to consider more detailed representations in the DES-OPF paradigm.

The following are noteworthy studies that may contribute further to the reader's understanding. Lv et al. [110] bring an interesting dimension to microgrid operation by considering the requirements of the distribution network within the optimisation framework. Their bi-level model has an upper level (master) where the distribution network minimises power losses and voltage offset through the output constraints of the microgrid and a lower level (slave) consisting of microgrids optimising operating cost. Their primary concern is that DES literature does not adequately consider the power limits at the Point of Common Coupling (PCC), i.e. where the microgrid connects to the distribution network. These limits are described below:

$$P^{min} \leq P_{PCC,t} \leq P^{max} \quad \forall t \qquad (17)$$

where $P^{min}$ and $P^{max}$ are pre-specified limits, and $P_{PCC,t}$ is the summation of the active powers at all nodes $n$ and time $t$. Van Leeuwen et al. [111] incorporate convexified branch flow equations in a blockchain-based operational model that considers trading aspects of a DES. While the study considers a balanced network and does not perform any detailed comparison with respect to DES-OPF, it demonstrates that nonlinear constraints can be incorporated into DES models that focus on peer-to-peer trading, artificial intelligence and blockchain applications.

*4.3  Summary of analysis*

To summarise, the optimisation approaches, simulation tools and type of formulation for power flow (nonlinear and nonconvex, convex relaxations, or linear approximations) are given in Table 1. The variety of DER technologies employed in these models are highlighted in

Table 2. The studies analysed under Section 4 show the importance of including detailed (and often nonlinear) representations of DES units and power flow. In the DES Design – OPF paradigm, all the papers analysed consider single-phase balanced networks, while in DES Operation, there is increasing interest in investigating multi-phase unbalanced networks. Hybrid optimisation approaches and/or simulation tools are commonly used in both design and operational studies, in an attempt to achieve a good accuracy-complexity balance. Some studies focus on obtaining globally optimal solutions while potentially compromising feasibility, while others focus on obtaining locally optimal but feasible solutions with respect to the nonlinear representations used. Regarding the use of tools, MATPOWER [92] has been used to simulate OPF in design studies that consider balanced networks, while OpenDSS [107] has been used in operational studies which consider three-phase unbalanced networks. Morstyn et al. [104] has introduced the *Network_3ph* class via the OPEN platform, which can be easily integrated to any Python-based DES operational models, either as an optimisation tool with multi-phase linearised formulations, or as a simulation tool with multi-phase nonlinear formulations. While OpenDSS remains an established simulation tool, comparisons between this and *Network_3ph* are promising. However, the applicability of *Network_3ph* to design models requires further investigation. Other tools offering nonlinear balanced and/or unbalanced power flow simulation and/or optimisation are summarised in Table 3. Note that there are many other tools available for OPF; however, they have not been included in this table as they only support linear power flow representations.

**Table 1.** The optimisation methods, simulation tools and power flow constraints employed in detailed and combined (DES-OPF) studies discussed.

| Ref. | Design | Operation | Optimisation approaches | Power flow Simulation tools | Balanced (B)/ Unbalanced Network | Type of formulation* |
|---|---|---|---|---|---|---|
| Morvaj et al. [8] | ✔ | - | GA, MILP | MATPOWER | B | N, L |
| Mashayekh et al. [78] | ✔ | - | MILP | - | B | L |
| Jordehi et al. [94] | ✔ | - | MINLP | - | B | N |
| Sfikas et al. [21] | ✔ | - | SQP | - | B | N |
| Nemati et al. [9] | - | ✔ | GA, MILP | MATPOWER | B | N, L |
| Chen et al. [99] | - | ✔ | MISOCP | - | B | C |
| Shuai et al. [12] | - | ✔ | ADP, DP, PSO | - | B | N |
| Zafarani et al. [101] | - | ✔ | NLP, LP | - | B | N, L |
| Shi et al. [11] | - | ✔ | PCPM | - | B | C |
| Thomas et al. [10] | - | ✔ | MILP | OpenDSS | U | L |
| Jin et al. [105] | - | ✔ | GA + MINLP | OpenDSS | U | N |
| Morstyn et al. [104] | - | ✔ | MILP, MINLP | PandaPower | B/U | N, L |
| Soares et al. [108] | - | ✔ | MILP + PSO | - | B | N, L |
| Lv et al. [109] | - | ✔ | GA + NLP | MATPOWER | B | N |
| van Leeuwen et al. [110] | - | ✔ | ADMM | - | B | C |

* Nonlinear and nonconvex formulations (N), Convex relaxations (C), and linear approximations (L)

Table 2. DER technologies used in the detailed models. The following notation is used: PV – Photovoltaic cells, CHP – Combined Heat and Power, GB – Gas Boilers, WT – Wind Turbines, DG – Diesel generators, FC – Fuel Cells, MT – Microturbines, ES – Electricity Storage, HS – Heat Storage, EN – Electricity Network, HN – Heat Network, EV – Electric Vehicles.

| Technology | PV | CHP | GB | WT | DG | FC | MT | ES | HS | EN | HN | EV | Other |
|---|---|---|---|---|---|---|---|---|---|---|---|---|---|
| Morvaj et al. [8] | ✔ | ✔ | ✔ | | | | | | ✔ | ✔ | | | |
| Mashayekh et al. [78] | ✔ | | | | | | ✔ | ✔ | | ✔ | ✔ | | Internal Combustion Engine, Chiller |
| Jordehi et al. [94] | ✔ | | | ✔ | | | | | | ✔ | | | Hydro, Geothermal, Battery Swap Station |
| Sfikas et al. [21] | ✔ | | | ✔ | | | | ✔ | | ✔ | | | Biomass |
| Nemati et al. [9] | ✔ | | | ✔ | ✔ | ✔ | ✔ | ✔ | | ✔ | | | |
| Chen et al. [99] | | ✔ | ✔ | ✔ | | | | ✔ | ✔ | ✔ | ✔ | | |
| Shuai et al. [12] | ✔ | | | ✔ | ✔ | | ✔ | ✔ | | ✔ | | | |
| Zafarani et al. [101] | | ✔ | | | | | | | | ✔ | ✔ | ✔ | |
| Shi et al. [11] | ✔ | | | ✔ | ✔ | | | ✔ | | ✔ | | | |
| Thomas et al. [10] | ✔ | | | ✔ | | | | ✔ | | ✔ | | ✔ | APS |
| Jin et al. [105] | ✔ | ✔ | ✔ | ✔ | | | | | | ✔ | ✔ | | Chiller |
| Morstyn et al. [104] | ✔ | | | | | | | ✔ | | ✔ | | ✔ | Flexible HVAC |
| Soares et al. [108] | ✔ | | | ✔ | | | | ✔ | | ✔ | | ✔ | Biomass |
| Lv et al. [109] | ✔ | ✔ | | ✔ | | | ✔ | ✔ | | ✔ | | | |
| van Leeuwen et al. [110] | ✔ | | | | | | | ✔ | | ✔ | | ✔ | |

Table 3. List of OPF tools with nonlinear model formulations.

| Tool | Simulation: Balanced OPF | Simulation: Unbalanced OPF | Optimisation: Balanced OPF | Optimisation: Unbalanced OPF | Notes |
|---|---|---|---|---|---|
| DiSC [112] | ✔ | - | - | - | Matlab-based |
| DIgSILENT PowerFactory [113] | ✔ | ✔ | ✔ | - | Not open-source |
| GridLab-D [114] | | ✔ | - | - | Can be linked to Matlab |
| MATPOWER [92] | ✔ | - | ✔ | - | Matlab-based |
| OpenDSS [107] | ✔ | ✔ | - | - | Can be linked to Matlab and Python |
| OATS [115] | ✔ | - | ✔ | - | Python-based |
| OPEN [104] | ✔ | ✔ | ✔ | ✔* | Python-based, uses PandaPower, unbalanced OPF through *Network_3ph* |
| PandaPower [109] | ✔ | ✔ | ✔ | - | Python-based |
| PyPSA [116] | ✔ | - | ✔* | - | Python-based |

* Denotes that the formulations are linearised.

## 5. Discussion and recommendations

Revisiting the three research questions highlighted in Section 1.1:

1. *What are the different levels of approximation evident through literature when modelling nonlinear aspects of a DES or microgrid?*

It is evident from Section 2 that DES baseline models predominantly use overly simplistic linear approximations, such as the DC approximations, and exclude detailed AC power flow constraints. A key reason for this is the additional complexity nonlinear constraints introduce, resulting in higher computational expense and no guarantee of global optimality. The distinct divide between OPF and DES baseline models can also be seen when examining the studies in Section 3. OPF models have higher granularity on network topology-related variables such as active and reactive power flows, voltages, and current, and simplified representations for DERs such as renewable generators and storage. On the other hand, DES models have greater detail on the design and operation of DERs and pay little attention to the underlying distribution network and associated constraints.

The detailed analysis presented in Section 4 affirms that OPF and DES models should not be considered as two entirely separate problems. Different levels of approximation are evident in these models as well. Most of the studies identified under the DES-OPF paradigm consider balanced networks, for which the representations are less complex and computationally expensive compared to unbalanced networks. However, in reality distribution networks are rarely balanced, therefore representation of unbalanced networks should, in theory, produce the most realistic and practically feasible results. Many of the studies in this section use linear approximations or convex relaxations of the nonlinear and nonconvex power flow constraints to achieve faster solution times and globally optimal solutions. However, there are insufficient comparisons between results of some of the approximations employed and those with original nonconvex formulations. This is evident in Design models, where different designs are reported between the two model formulations, but the potential reasons for these differences are not investigated. Furthermore, there is greater focus on achieving globally optimal solutions; however, whether these solutions are practically feasible has not been investigated thoroughly. One way of confirming the feasibility of a design is to test it using a high-fidelity simulation tool, such as those recommended in Table 3. This will ensure that the designs are suitable for implementation, where relevant.

2. *Should detailed power flow constraints be included in DES models?*

The detailed and combined models in Section 4 highlight that it is essential to include detailed power flow constraints in DES design and operational models that are connected to the AC grid, such as those mentioned in Table 1. Some of these studies demonstrate that the inclusion of AC OPF constraints can have significant impact on the objective functions and resulting designs when compared to models that employ less detailed representations.

3. *Where applicable, how are these different levels linked within models to produce optimal solutions within a reasonable time frame?*

The detailed and combined optimisation models analysed in Section 4 identify simulation tools which can be linked to optimisation models to provide an additional layer of detail, especially when considering nonlinear and nonconvex power flow constraints. The capabilities of these tools are summarised in Table 3. It typically involves using a linear approximation within the optimisation model, combined with the use of the simulation tool to verify the calculations or detect violations at each iteration. This is because many of these tools have not been adapted for optimisation models. Furthermore, as indicated in Table 3, while there are many simulation tools available, only a few of these offer multi-phase, unbalanced, and nonlinear OPF.

Although a wide range of metaheuristic techniques have been applied in recent times, and especially in DES-OPF studies, it is possible and preferable to use deterministic mathematical programming approaches to solve detailed models. The latter continue to provide better and more reliable solutions compared to metaheuristics.

With respect to the more subjective research question posed:

*Out of the modelling approaches explored, what level or method is most suitable to achieve a good complexity-accuracy balance for DES - OPF?*

we have identified that the following considerations are necessary when integrating power flow in DES design and operation:

1. Does the model consider a balanced or an unbalanced network?
2. Which power flow formulation best represents the balanced/unbalanced network (bus injection, branch flow, Z-bus, etc.)?
3. What additional or "side" constraints, such as current constraints, transformer constraints, and thermal limits, will be included in the OPF formulation?
4. What are the nonlinearities associated with DES technologies that can further increase complexity of the model and its solution?
5. Are there discrete variables in the model, and are there any that can be eliminated without compromising the accuracy of the DER representation?
6. If the priority is to obtain globally optimal solutions or obtain fast solutions, which linear approximations or convex relations can be used?
7. If appropriate, what high-fidelity tools can be used to verify the results from approximations?

These discussion points lead us to recommend that future studies should aim to consolidate DES design and OPF, and investigate multi-phase unbalanced networks. Furthermore, comparisons are required between high-fidelity models and models with linear approximations to understand what level of detail is required for day-ahead and real-time operation of DES. Finally, simulation tools, such as those in Table 3, are a good starting point for DES modellers who have little knowledge on power flow studies and wish to incorporate detailed network constraints.

## 6. Conclusions

The integration of DES and microgrids into existing energy networks has become increasingly important, as many countries aim to reduce carbon emissions and meet increasing demand through renewable energy resources. Although there are many studies investigating the design and operation of DES, especially using optimisation models, they often either ignore or linearise inherent nonlinearities related to generation/storage technologies, and the underlying distribution network for power flow. Such models produce results that are often sub-optimal or even practically infeasible, meaning that they cannot be implemented. Furthermore, OPF studies that include detailed nonlinear power flow constraints are considered as standalone problems, thus leading to the disparity between OPF and DES design and operational models.

We identify a new class of models, labelled DES-OPF models, which bridge the gap between DES design and operation and OPF by considering relevant nonlinearities and details in both aspects. By identifying the characteristics of DES baseline models and OPF base formulations first, the review sheds light on technical constraints related to OPF which should be included in DES design and operational models. These have often been overlooked in the latter and replaced with oversimplified linear approximations such as the DC approximations. Analysis of DES-OPF models shows that the use of innovative optimisation approaches and integration of simulation tools can undoubtedly help achieve a better accuracy-complexity balance, without excluding vital nonlinear constraints. Furthermore, to guide future efforts towards the consolidation of DES and OPF approaches, our considerations and recommendations stress the importance of conducting further research in this area and highlight the need for more comparisons between high-fidelity models and models with linear approximations.

To conclude, there is increasing focus on understanding the multiple facets of DES (not limited to optimal design, operation, and power flow), and producing models that adequately represent each aspect. This forms the basis of the accuracy-complexity balance modellers seek and can be individual to each DES or microgrid. Gaining better understanding of each of the facets and attempting to consolidate them will undoubtedly lead to the development of more robust and representative models that can be implemented in reality.